\documentclass[fleqn,a4paper,12pt]{article} 

\usepackage{amsmath}
\usepackage{amssymb}
\usepackage{rotate}
\newtheorem{theorem}{Theorem}
\newtheorem{proposition}[theorem]{Proposition}
\newtheorem{corollary}[theorem]{Corollary}

\begin{document}

\title{$DW(2n,q)$, $n\geq 3$, has no ovoid:\\
A single proof}
\author{Harm Pralle}
\maketitle

\begin{abstract}
An ovoid of a dual polar space $\Delta$ is a point set meeting every
line of $\Delta$ in exactly one point. For the symplectic dual polar
space $DW(6,q)$, Cooperstein and Pasini \cite{CP2003} have recently proved
no ovoid exists if $q$ is odd. Earlier, Shult has proved
the same for even $q$ (cf.\ \cite[2.8]{PS2001}). In this paper, we
prove the non-existence of ovoids in $DW(6,q)$ independently from the
parity of $q$. 
\end{abstract}

{\footnotesize {\bf MSC 2000:} 51A15, 51A50, {\bf Key words:} dual polar
spaces, hyperplanes, ovoids.}

\bigskip
\section{Introduction}
Let $\Delta$ be a finite dual polar space of finite rank $n\geq 3$, i.e.\ it is
the geometry dual of a polar space $\Pi$. The points of $\Delta$ are
the $(n-1)$-dimensional singular subspaces of $\Pi$, the lines of
$\Delta$ are the $(n-2)$-dimensional singular subspaces of $\Pi$ and,
more generally, the elements of type $i$ of $\Delta$, $i=1,...,n$, are
the $(n-i)$-dimensional singular subspaces of $\Pi$. Thus $\Delta$
belongs to the diagram

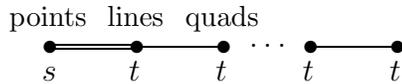
\parindent 2em
\setlength{\unitlength}{0.07em}
\begin{picture}(110,43)
\linethickness{0.04em}
\put(15,12){$\bullet$}
\put(18,14.5){\line(1,0){38}}
\put(18,16.5){\line(1,0){38}}
\put(55,12){$\bullet$}
\put(58,15.5){\line(1,0){40}}
\put(95,12){$\bullet$}
\put(108,15.5){${\bf \ldots}$}
\put(135,12){$\bullet$}
\put(138,15.5){\line(1,0){40}}
\put(175,12){$\bullet$}
\put(15,0){$s$}
\put(55,0){$t$}
\put(95,0){$t$}
\put(135,0){$t$}
\put(175,0){$t$}
\put(0,25){{\small points}}
\put(45,25){{\small lines}}
\put(81,25){{\small quads}}
\end{picture}
\parindent 0em

An ovoid of a dual polar space is a point set that meets each line in
exactly one point. It is an outstanding conjecture that in
finite classical dual polar spaces of rank at least three no ovoid
exists. By Pasini and Shpectorov \cite{PS2001}, the complement of an
ovoid of a dual polar space of rank $3$ cannot be
flag-transitive. Recently, Cooperstein and Pasini \cite{CP2003} have
proved that the finite symplectic dual polar space $DW(6,q)$ for $q$
odd has no ovoid. Together with a similar non-existence result for $q$
even by Shult to be found in \cite[2.8]{PS2001}, $DW(6,q)$ has no ovoid. 

\bigskip
Since we restrict the considerations to finite dual polar spaces
$\Delta$, the polar space $\Pi$ dual of $\Delta$ is classical by Tits'
classification of polar spaces. Suppose $O$ is an ovoid of
$\Delta$. The ovoid $O$ intersects each quad of $\Delta$ in an ovoid
of a generalized quadrangle. The quads are either grids if $t=1$
or classical generalized quadrangles of order $(s,t)$. Since we assume they
admit ovoids and lines aren't short, i.e.\ $t\geq 2$, $\Delta$ is the
dual of the symplectic polar space $W(6,q)$, the orthogonal polar
space $O^-(8,q)$ or possibly the hermitian polar space $H(7,q^2)$,
$q\geq 3$, with the quads being generalized quadrangles $O(4,q)$,
$H(4,q^2)$ or the dual of $H(5,q^2)$, respectively (for the classical
generalized quadrangles admitting ovoids, see Payne and Thas
\cite{PT1984}). Note that the existence of an ovoid of the dual of the
hermitian generalized quadrangle $H(5,q^2)$, $q\geq 3$, i.e.\ a spread
of $H(5,q^2)$, is an open problem.

\medskip
In this paper, we simplify the proof of \cite{CP2003} and generalize
it such that we can apply it to any dual polar space of rank $3$
hypothetically admitting ovoids. For the classical dual polar spaces,
the counting leads to a contradiction only for the symplectic dual
polar space $DW(6,q)$. However, it is independent from the parity of
$q$, whence our new prove comprises the results of Shult
(\cite[2.8]{PS2001}) and Cooperstein and Pasini \cite{CP2003}.

\begin{theorem} \label{Thm}
If $\Delta$ has an ovoid, then $s^2-s-t^2-t \geq 0$. In particular,
$DW(6,q)$ has no ovoid.
\end{theorem}

Theorem \ref{Thm} generalizes to symplectic dual polar spaces
$DW(2n,q)$ of finite rank $n\geq 4$ since the point-line residue of
each element of type $4$ of $DW(2n,q)$ is a dual polar space $DW(6,q)$
which would intersect an ovoid of $DW(2n,q)$ in an ovoid of $DW(6,q)$
in contradiction to Theorem \ref{Thm}.

\begin{corollary}
A finite symplectic dual polar space $DW(2n,q)$, $n\geq 3$, has no
ovoid. \hfill $\Box $
\end{corollary}

\section{Proof of the Theorem}

We recall the notation and introduce some more terminology. For an
element $x$ of a geometry $G$ of diameter $d$, for $i=1,...,d$,
$G_i(x)$ denotes the set of points of $G$ at distance $i$ from $x$ in the
collinearity graph of $G$.

\smallskip
Before starting the proof, we introduce the projection $\pi_E$ of the
point set of a dual polar space onto the point set of an arbitrary
element $E$ of type at least $2$. Since dual polar spaces of rank
$n$ are near $2n$-gons (cf.\ Cameron \cite{C1982}), for any point $p$
and any element $E$ of the dual polar space not on $p$, there is
a unique point $\pi_E(p)$ in $E$ nearest $p$. In setting $\pi_E(x)=x$
for all points $x$ in $E$, the so-defined mapping $\pi_E$ maps two
collinear points of $\Delta$ either onto the same point of $E$ or onto two
collinear points of $E$. Thus, if $\Delta$ is a dual polar space of
rank $3$, the following proposition which mainly serves as reference
for Proposition \ref{mu}, is immediate.

\begin{proposition} \label{projection}
If $\omega$ is a quad of $\Delta$, then a line $l$
disjoint from $\omega$ is mapped onto the line $\pi_\omega(l)$ of
$\omega$ such that for each point $p$ on $l$, there is a unique point
on $\pi_\omega(l)$ collinear with $p$. 
\hfill $\Box $
\end{proposition}

Let us now turn to the proof of Theorem \ref{Thm}. We suppose
$\Delta$ is a finite dual polar space of rank $3$ and the generalized
quadrangle consisting of the points and lines of a quad of $\Delta$
has order $(s,t)$ with $t\geq 2$. Moreover, since we assume $\Delta$
has an ovoid $O$ and the quads are classical generalized quadrangles,
it follows $s\geq t$ (cf.\ Payne and Thas \cite{PT1984}).

\medskip
Let $\infty$ be a point of $O$. Then the set $H$ of points of $\Delta$
at non-maximal distance from $\infty$ is the so-called {\em singular
hyperplane with deepest point} $\infty$, whence 
%. E.g.\ the singular hyperplane with deepest point $\infty$ is the set 
$H = \{ \infty \} \cup \Delta_1(\infty) \cup \Delta_2(\infty)$. Denote
the affine dual polar space consisting of the elements of $\Delta $
not contained in $H$ by $\Gamma := \Delta - H$, i.e.\ the point set of
$\Gamma$ is $\Delta_3(\infty)$. We call the lines of $\Delta$ not
contained in the singular hyperplane $H$ the {\em affine lines} of
$\Gamma$. Moreover, we set $\Omega := O\cap \Gamma$.

Since both $O$ and $H$ are hyperplanes of $\Delta$, a line $l$ of
$\Delta$ not contained in $H$ has exactly one point $l^O = l\cap O$ of
the ovoid $O$ and one point $l^\infty = l\cap H$ of $H$. 
The main idea of Cooperstein and Pasini \cite{CP2003} is to count
pairs $(l,m)$ of concurrent affine lines of the affine dual polar
space $\Gamma = \Delta - H$ such that the unique point $l^\infty$ of
the line $l$ of $\Delta$ in $H$ lies in $O$, i.e.\ $l^\infty = l^O$,
and the unique point $m^\infty$ of the line $m$ of $\Delta$ in $H$
does not belong to $O$, i.e.\ $m^\infty \neq m^O$. Using Cauchy's
inequality, their final conclusion for $\Delta = DW(6,q)$, $q$ odd, is
$2q \leq 0$ proving no ovoid exists.

\smallskip
We follow most of the proof of Cooperstein and Pasini \cite{CP2003}
and most of their notation. The only modification is the method to
prove  Proposition \ref{mu} below. It allows us to
apply the final argument of \cite{CP2003} to arbitrary
finite dual polar spaces. In particular, it includes $DW(6,q)$ for $q$
even.

\begin{proposition}
It holds $|\Gamma|=s^3t^3$, $|O|=(st+1)(st^2+1)$, and $|\Omega|=st^3(s-1)$.
\end{proposition}

{\bf Proof.}\quad Since $|\Delta| = (t+1)(st+1)(st^2+1)$,
$|\Delta_1(\infty)| = (t^2+t+1)s$ and
$|\Delta_2(\infty)|=|\Delta_1(\infty)|(t^2+t)s/(t+1) = (t^2+t+1)st$,
it follows 
\[ |\Gamma|=|\Delta|-|\Delta_2(\infty)|-|\Delta_1(\infty)|-1
= s^3t^3 \ .\] 
We determine $|O|$ by
\[ |O|=\frac{\# \mbox{ lines}%\cdot \# \mbox{ ovoid pointsper line}
}{\# \mbox{ lines per point}}= \frac{(st+1)(st^2+1)(t^2+t+1)}{t^2+t+1} = (st+1)(st^2+1) \ .\]
Since each quad on $\infty$ has $st$ points of $O-\{\infty\}$ and none
of the points of $O$ at distance two from $\infty$ belongs to two
quads on $\infty$, it follows 
\[ |\Omega| = |O| - (t^2+t+1)st -1= st^3(s-1) \ . 
\hfill \Box \]

\bigskip
Following Cooperstein and Pasini \cite{CP2003}, for a point $p\in
{\cal G} := \Gamma - \Omega$ let
\[ \mu_p := |\Gamma_1(p) \cap \Omega | 
\]
be the number of ovoid points collinear with $p$ not belonging to the singular
hyperplane $H$. Since there are $t^2+t+1$ lines on $p$ all
meeting $O$ in exactly one point, $t^2+t+1-\mu_p$ lines on $p$ meet
$O$ in points of $O-\Omega$. Then the number of pairs of concurrent
lines of $\Gamma$ one meeting $\Omega$ and the other meeting $O- \Omega$ is
\[ N := \sum_{p \in {\cal G}} \mu_p (t^2+t+1-\mu_p ) 
\]
since any two such lines meet in a point of ${\cal G}$. 

We determine $N$ in the following. Denote by $\cal L$ the set of
affine lines of $\Gamma$ not meeting $\Omega$ and by $\cal M$ the set
of affine lines of $\Gamma$ meeting $\Omega$. Let $l\in {\cal L}$ and
$m\in {\cal M}$, i.e.\ $l^\infty \in O\cap H = O-\Omega$ and $m^\infty
\in H-O$. We set
\begin{eqnarray*} 
\mu^-(l) &:=& \sum_{p\in l} \mu_p \mbox{\ \ and }\\
\mu^+(m) &:=& \sum_{p\in m\cap {\cal G}} \mu_p \ .
\end{eqnarray*}
The number $\mu^-(l)$ is the number of affine lines of $\Gamma = \Delta
-H$ concurrent with $l$ and meeting $O$ in a point of $\Omega = \Gamma \cap
O$. Then the number $N$ of pairs $(l,m)$ of concurrent affine
lines with $l\in {\cal L}$ and $m\in {\cal M}$ is
\[N=\sum_{l\in {\cal L}} \mu^-(l) \ .\] 
The following proposition determines the numbers $\mu^-(l)$ and
$\mu^+(m)$ showing the numbers are independent from the choice of the
particular lines $l\in {\cal L}$ and $m\in {\cal M}$. Note that the
number $\mu^+(m)$ will be used lateron, too.

\begin{proposition} \label{mu}
For $l\in {\cal L}$ and $m \in {\cal M}$, it holds
$\mu^-(l) = (s-1)(t^2+t)$ and $\mu^+(m) = (s-1)(t^2+t+1)-(t^2+t)$.
\end{proposition}

{\bf Proof.}\quad To count the points of $\Omega$ collinear with $l$,
we have to subtract the number $M$ of points of $O\cap H$ collinear
with $l$ from the number $s(t^2+t)$ of all points of $O$ collinear
with $l$ in $\Gamma$. Since in $\Delta$, $l$ meets $\Delta_2(\infty)$
in the point $l^\infty = l\cap H$, there is a unique quad $\delta$ on
$\infty$ meeting $l$, namely in $l^\infty$. For each quad $\omega \neq
\delta$ on $\infty$, the affine line $l$ is projected by
$\pi_{\omega}$ onto the line $\pi_{\omega}(l)$ which does not go
through $\infty$ (cf.\ Proposition \ref{projection}). Hence
$\pi_{\omega}(l)$ contains exactly one point of $(O\cap H) -
\{\infty\}$. For any two quads $\omega, \tau \neq \delta $ on
$\infty$, the points $\pi_{\omega}(l) \cap O$ and $\pi_{\tau}(l)\cap
O$ are distinct. Hence there are $t^2+t$ points of $O\cap H$ collinear
with $l$. It follows  
\[ \sum_{p\in l} \mu_p = s(t^2+t) - (t^2+t)= (s-1)(t^2 + t) \ . \]
Similarly, let $\gamma$ be the quad on $\infty$ meeting the line
$m$. As before, for each quad $\sigma \neq \gamma$ on $\infty$,
$\pi_{\sigma}(m)$ contains a unique point of $(O\cap
H)-\{\infty\}$. Since the affine part of $m$ has only $s-1$ points of
$\cal G$ and there are $t^2+t$ quads on $\infty$ distinct from
$\gamma$, it follows 
\[ \sum_{p\in m\cap {\cal G}} \mu_p = (s-1)(t^2+t+1) - (t^2+t) 
\hfill \Box \]

\bigskip
In particular, $\mu^-(l)$ and $\mu^+(m)$ do not depend on the lines
$l$ and $m$. Thus $N=\sum_{l\in {\cal L}} \mu^-(l) = |{\cal L}| \cdot
\mu^-(l)$ for some line $l\in {\cal L}$. Since $\cal L$ is the set of
affine lines meeting $H\cap O$ and since $|H\cap O| = st(t^2+t+1)$, it
follows $|{\cal L}|=st^3(t^2+t+1)$. Thus the number $N$ follows:

\begin{corollary} \label{Cor N}
$N=st^4(s-1)(t+1)(t^2+t+1)$
\hfill $\Box$
\end{corollary}

To determine $\sum_{p\in {\cal G}} \mu_p$, consider the following partition of
the point set of ${\cal G}=\Gamma-O$. Let $\kappa$ be a quad on $\infty$. Then
the affine lines of $\Gamma$ meeting $\kappa$, considered as lines of
$\Delta$, partition the point set of $\Gamma$. These affine lines fall
in two classes $K^-$ and $K^+$ where $K^-$ is the set of affine lines
meeting $\kappa$ in points of $O - \{\infty\}$, i.e.\ their
affine part does not meet $\Omega$, and $K^+$ is the set
of affine lines meeting $\kappa$ in points of $(\kappa
-\infty^{\perp}) - O$, i.e.\ their affine part meets $\Omega$. Since
$|\kappa \cap O|=st+1$, it follows $|(\kappa -\infty^{\perp}) -
O| = st(s-1)$. Thus $|K^-|=st^3$ and $|K^+|=st^3(s-1)$. It follows
\begin{eqnarray*} 
\sum_{p\in {\cal G}} \mu_p &=& |K^-| \mu^- + |K^+|\mu^+ \\
&=& st^3(s-1)(t^2+t)+st^3(s-1)((s-1)(t^2+t+1)-(t^2+t))\\
&=& st^3(s-1)^2(t^2+t+1)
\end{eqnarray*}
With Corollary \ref{Cor N} it follows
\begin{eqnarray*} 
\sum_{p\in {\cal G}} \mu_p^2 &=& (t^2+t+1) \sum_{p\in {\cal G}} \mu_p - N \\
&=& (t^2+t+1)st^3(s-1)^2(t^2+t+1) - st^4(s-1)(t+1)(t^2+t+1) \\
&=& st^3(t^2+t+1)(s-1) ((t^2+t+1)(s-1)-(t^2+t))
\end{eqnarray*}
We are now in the position to conclude the proof by Cauchy's inequality
\[ |{\cal G}| \sum_{p\in {\cal G}} \mu_p^2 \geq \left( \sum_{p\in {\cal G}}
\mu_p \right)^2 \ . 
\]
Indeed, it holds
\[ |{\cal G}| \sum_{p\in {\cal G}} \mu_p^2 =
s^2t^6 (s^2-s+1) (t^2+t+1)(s-1) ((t^2+t+1)(s-1)-(t^2+t))
\]
leading with $\sum_{p\in {\cal G}} \mu_p = st^3(s-1)^2(t^2+t+1)$ to
the inequality 
\[(s^2-s+1)((t^2+t+1)(s-1)-(t^2+t)) \geq (s-1)^3(t^2+t+1) \ .\]
This is equivalent to
\[ s^2-s-t^2-t \geq 0 \ .\]
Note that the counting arguments do not depend on the particular
polar space under consideration. In particular, they are independent
from the parity of $s$ or $t$.

Since the existence of an ovoid in a classical generalized quadrangle forces
$s\geq t$, %(cf.\ Payne and Thas \cite{PT1984})
the inequality leads to a contradiction only if $s=t$. From the dual
polar spaces admitting ovoids in quads, only $\Delta \cong DW(6,q)$
has orders $s=t=q$ leading to the contradiction $2q\leq 0$. Hence the
dual polar space $DW(6,q)$ does not admit any ovoid whereas we cannot
deduce anything about the two other classical dual polar spaces
$DO^-(8,q)$ admitting ovoids of quads and $DH(7,q^2)$ possibly
admitting ovoids of quads.

\bigskip
H.\ Pralle\\
Institute for Geometry, Algebra and Discrete Mathematics\\
TU Braunschweig\\
Pockelsstr.\ 14\\
38106 Braunschweig, Germany\\
H.Pralle@tu-bs.de

\bigskip
MSC 2000: 51A15, 51A50

\bigskip
Key words: dual polar spaces, hyperplanes, near $2n$-gons, ovoids.

\pagebreak


\begin{thebibliography}{1}

\bibitem{C1982}
Cameron, P.~J.
\newblock Dual Polar Spaces.
\newblock {\em Geom.~Dedicata} {\bf 12} (1982), 75--85.

\bibitem{CP2003}
Cooperstein, B., and Pasini, A.
\newblock $DW(6,q)$ has no ovoid.
\newblock {\em J.\ Combin.\ Th.\ A}, to appear.

\bibitem{PS2001}
Pasini, A., and Shpectorov, S.
\newblock Uniform hyperplanes of finite dual polar spaces of rank 3.
\newblock {\em J.~Comb.~Theory (A)} {\bf 94} (2001), 276--288.

\bibitem{PT1984}
Payne, S., and Thas, J.~A.
\newblock {\em Finite Generalized Quadrangles}.
\newblock Pitman, Boston (1984).

\end{thebibliography}
\end{document}